\documentclass[A4paper,11pt,reqno]{amsart}

\usepackage[english]{babel}
\usepackage{amsmath}
\usepackage{amssymb}
\usepackage{amsfonts}
\usepackage[all,color]{xy}
\usepackage{hyperref}
\usepackage{lscape}
\usepackage{graphicx}
\usepackage{float}

\numberwithin{equation}{section}
\usepackage{color}
\usepackage{hyperref}
\usepackage{enumerate,xspace}
\usepackage{graphicx}

\definecolor{lightgray}{rgb}{0.666666,0.666666,0.666666}
\definecolor{pinegreen}{rgb}{0.15,0.7,0.15}
\definecolor{darkgreen}{rgb}{0,0.45,0}

 \definecolor{lightgrey}{rgb}{0.666666,0.666666,0.666666}

\renewcommand{\phi}{\varphi}
\renewcommand{\epsilon}{\varepsilon}

 \newtheorem{proposition}{Proposition}[section]
 
 \newtheorem*{lemma*}{Lemma}
 \newtheorem{theorem}[proposition]{Theorem}
 
 \theoremstyle{definition}
 \newtheorem{construction}[proposition]{Construction}
 
 \newtheorem*{definition*}{Definition}
 \newtheorem{example}[proposition]{Example}
 \newtheorem*{example*}{Example}
 \newtheorem{remark}[proposition]{Remark}

\numberwithin{equation}{section}

\DeclareMathOperator{\Mnd}{Mnd}

\DeclareMathOperator{\Decr}{Dec_r}
\DeclareMathOperator{\Decl}{Dec_l}
\DeclareMathOperator{\ev}{ev}

\DeclareMathOperator{\End}{End}
\DeclareMathOperator{\BS}{BS}

\newcommand{\dtwocell}[3][0.5]{\ar@{}[#2] \ar@{=>}?(#1)+/u  0.2cm/;?(#1)+/d 0.2cm/^{#3}}

\newcommand{\cc}{\ensuremath{\mathcal C}\xspace}

\newcommand{\ck}{\ensuremath{\mathcal K}\xspace}

\newcommand{\cv}{\ensuremath{\mathcal V}\xspace}

\newcommand{\bbm}{\ensuremath{\mathbb M}\xspace}
\newcommand{\bbn}{\ensuremath{\mathbb N}\xspace}

\newcommand{\bbz}{\ensuremath{\mathbb Z}\xspace}

\newcommand{\bfk}{\ensuremath{\mathbf K}\xspace}

\newcommand{\ox}{\otimes}
\newcommand{\x}{\times}

\newcommand{\op}{\ensuremath{{}^{\textrm{op}}}\xspace}
\newcommand{\rev}{\ensuremath{{}^{\textrm{rev}}}\xspace}

\newcommand{\lMod}{\ensuremath{\textnormal{-}\mathbf{Mod}}\xspace}

\newcommand{\bMod}{\ensuremath{\textnormal{-}\mathbf{Mod}\textnormal{-}}\xspace}

\newcommand{\NLax}{\ensuremath{\mathbf{NLax}}\xspace}

\newcommand{\Set}{\ensuremath{\mathbf{Set}}\xspace}
\newcommand{\SSet}{\ensuremath{[\Delta\op,\Set]}\xspace}
\newcommand{\Cat}{\ensuremath{\mathbf{Cat}}\xspace}

\newcommand{\Ab}{\ensuremath{\mathbf{Ab}}\xspace}

\begin{document}

\title{Hochschild homology, lax codescent,  and duplicial structure}

 \author{Richard Garner}
\address{Department of Mathematics, Macquarie University NSW 2109, Australia}
\email{richard.garner@mq.edu.au}

\author{Stephen Lack}
\address{Department of Mathematics, Macquarie University NSW 2109, Australia}
\email{steve.lack@mq.edu.au}

\author{Paul Slevin}
\address{University of Glasgow, School of Mathematics and Statistics, 15 University Gardens, G12 8QW Glasgow, UK}
\email{p.slevin.1@research.gla.ac.uk}

\begin{abstract}
We study the duplicial objects of Dwyer and Kan, which generalize the cyclic objects of Connes. We describe duplicial objects in terms of the decalage comonads, and we give a conceptual account of the construction of duplicial objects due to B\"ohm and \c Stefan.
This is done in terms of a 2-categorical generalization of Hochschild homology. We also study duplicial structure on nerves of categories, bicategories, and monoidal categories.
\end{abstract}
\date\today
\maketitle

\section{Introduction}

The cyclic category $\Lambda$ was introduced by Connes \cite{Connes-cyclic} as part of his program to study non-commutative geometry. Cyclic objects, given by functors with domain $\Lambda$, have been studied by too many authors to list here, but many of these can be found in the reference list of the classic book \cite{Loday-cyclic} of Loday.

Various generalizations of cyclic structure have been considered; in particular the notion of {\em duplicial object} was studied by Dwyer and Kan in \cite{DwyerKan-cyclic}. These are given by functors with domain $\bfk\op$, for a certain category $\bfk$ of which $\Lambda$ is a quotient. Like cyclic objects, duplicial objects are simplicial objects equipped with extra structure. In both cases, the extra structure involves an endomorphism $t_n\colon X_n\to X_n$ of the object of $n$-simplices, for each $n$, subject to various conditions relating it to the simplicial structure. The difference between the two notions is that in the case of cyclic structure, the map $t_n$ is an automorphism of order $n+1$, so that $t^{n+1}_n=1$.

There is also an intermediate notion, in which the $t_n$ are required to be invertible, but the condition that $t^{n+1}_n=1$ is dropped. This was called {\em paracyclic structure} in \cite{GetzlerJones-paracyclic}, and also studied in \cite{Elmendorf-CyclicDuality} where the indexing category was called the ``linear category''. Somewhat confusingly, the name paracyclic has also been used by some authors to refer to what is called duplicial by Dwyer and Kan.

In this paper we shall provide a new perspective on duplicial structure, as well as analyzing ways in which it arises.

As explained, for example, in \cite{CWM}, a comonad on a category gives rise to simplicial structure on each object of that category, and this is the starting point for many homology theories. Just as simplicial structure can be used to define homology, cyclic (or duplicial or paracyclic) structure  can be used to define cyclic homology. In a series of papers \cite{BohmStefan1,BohmStefan2,BohmStefan3}, B\"ohm and \c Stefan looked at what further structure than a comonad is needed to equip the induced simplicial object with duplicial structure; the main extra ingredient turned out to be a second comonad with a distributive law \cite{Beck:dist-law} between the two. They also showed that their machinery could be used to construct the cyclic homology of bialgebroids. This was further studied in the papers \cite{Slevin1,Slevin2} by the third of us, along with various coauthors.

In the case of comonads and simplicial structure, there is a universal nature to the construction, once again explained in \cite{CWM}, and also in Section~\ref{sect:simplicial} below. There is no analogue given in the analysis of B\"ohm-\c Stefan, and our first goal is to provide one.

As well as the construction of simplicial structure from comonads, we also consider a second way that simplicial structure arises, namely as nerves of categories or other (possibly higher) categorical structures. Our second main goal is to analyze when the simplicial sets arising as nerves can be given duplicial structure.

The third main achievement of the paper actually arose as a by-product of our investigations towards the first goal. It is a connection between duplicial structure, especially as arising via the B\"ohm-\c Stefan construction, and Hochschild homology and cohomology. We shall present this first. We consider some very simple aspects of Hochschild homology and cohomology, only involving the zeroth homology and cohomology, and we generalize it to a 2-categorical context in a ``lax'' way. The resulting theory allows us to recapture the B\"ohm-\c Stefan construction as a sort of cap product in a very special case.

\subsection*{Acknowledgements}
The first two authors acknowledge support of an Australian Research Council Discovery Grant DP130101969; Garner was also supported by an Australian Research Fellowship DP110102360 and Lack by a Future Fellowship FT110100385. Slevin acknowledges an EPSRC Doctoral Training Award and funding from the University of Glasgow and Macquarie University which enabled his visit to Sydney in April--July 2015, during which period the bulk of this research was carried out. We are grateful to Tony Elmendorf for helpful discussions. 

\section{Simplicial structure, comonads, and decalage}
\label{sect:simplicial}

In this section we recall various ideas related to simplicial structure, most of which are well-known, although the notation used varies. The one new result is Proposition~\ref{prop:duplicial-decalage}, which reformulates the notion of duplicial structure in terms of the decalage comonads.

\subsection{Simplicial structure arising from comonads}\label{sect:bar}

We write \bbm for the strict monoidal category of finite ordinals and order-preserving maps, with tensor product given by ordinal sum and the empty ordinal serving as the unit. This is sometimes known as the ``algebraists' $\Delta$'', and is denoted by $\Delta$ in \cite{CWM} and $\Delta_+$ in many other sources, such as \cite{Verity-complicial}.

The full subcategory of \bbm consisting of the non-empty finite ordinals is isomorphic to the usual $\Delta$ (the ``topologists' $\Delta$''). A contravariant functor defined on $\Delta$ is a simplicial object, while a contravariant functor defined on (the underlying category of) \bbm is an augmented simplicial object.

\bbm is the ``universal monoidal category containing a monoid'', in the sense that for any strict monoidal category \cc, there is a bijection between monoids in \cc and strict monoidal functors from \bbm to \cc. (Similarly, if \cc is a general monoidal category then to give a monoid in \cc is equivalent, in a suitable sense, to giving a strong monoidal functor from \bbm to \cc.)

Dually, there is a bijection between comonoids in \cc and strict monoidal functors from $\bbm\op$ to \cc, and so any comonoid in \cc determines an augmented simplicial object in \cc. In particular, we could take \cc to be the strict monoidal category $[X,X]$ of endofunctors of a category $X$, so that a comonoid in \cc is just a comonad on $X$. Then any comonad $g$ on $X$ determines a unique strict monoidal functor $\bbm\op\to[X,X]$. We may now transpose this so as to obtain a functor $X\to [\bbm\op,X]$ sending each object of $X$ to an augmented simplicial object in $X$ called its \emph{bar resolution} with respect to $g$.

When, in the introduction, we referred to the ``universal nature'' of the construction of simplicial objects from comonads, it was precisely this analysis, using the universal property of $\bbm$, which we had in mind, and which we shall extend so as to explain the B\"ohm-\c Stefan construction.

\begin{remark}\label{rmk:Delta-rev}
There is an automorphism of $\bbm$ which arises from the fact that the opposite of the ordinal
$$n=\{0<\ldots<n-1\}$$
is isomorphic to $n$ itself. The automorphism fixes the objects, and sends an order-preserving map $f\colon m\to n$ to $f\rev$, where $f\rev(i)=m-1-f(n-1-i)$. This automorphism reverses the monoidal structure, in the sense that $n+n'=n'+n$ on objects, while for morphisms $f\colon m\to n$ and $f'\colon m'\to n'$ we have $(f+f')\rev=(f')\rev+f\rev$.
\end{remark}

\subsection{The decalage comonads}

The monoidal structure on \bbm extends, via Day convolution \cite{Day-convolution}, to a monoidal structure on the category $[\bbm\op,\Set]$ of augmented simplicial sets. The resulting structure is non-symmetric, but closed on both sides, so that there is both a left and a right internal hom.

Since the ordinal $1$ is a monoid in \bbm, the representable $\bbm(-,1)$ is a monoid in $[\bbm\op,\Set]$, and so on taking the internal hom out of $\bbm(-,1)$ becomes a comonad; or rather, there are two such comonads depending on whether one uses the left or right internal hom. These are called the decalage comonads, and they both restrict to give comonads, also called decalage, on the category $[\Delta\op,\Set]$ of simplicial sets.

As well as this abstract description, there is also a straightforward explicit description, which we now give for the case of augmented simplicial sets.

Given an augmented simplicial set $X$ as in the diagram
$$\xymatrix{
\llap{\ldots~~} X_2 \ar@<3ex>[r]|{d_0} \ar[r]|{d_1} \ar@<-3ex>[r]|{d_2} &
X_1 \ar@<1.5ex>[l]|{s_1} \ar@<-1.5ex>[l]|{s_0}  \ar@<1.5ex>[r]|{d_0} \ar@<-1.5ex>[r]|{d_1} &  X_0 \ar[l]|{s_0} \ar[r]^-{d_0} & X_{-1}}$$
the right decalage $\Decr(X)$ of $X$ is the augmented simplicial set
$$\xymatrix{
\llap{\ldots~~}  X_3 \ar@<3ex>[r]|{d_0} \ar[r]|{d_1} \ar@<-3ex>[r]|{d_2} &
X_2 \ar@<1.5ex>[l]|{s_1} \ar@<-1.5ex>[l]|{s_0}  \ar@<1.5ex>[r]|{d_0} \ar@<-1.5ex>[r]|{d_1} &  X_1 \ar[l]|{s_0} \ar[r]^-{d_0} & X_0}$$
obtained by discarding $X_{-1}$ and the last face and degeneracy map in each degree. There is a canonical map $\epsilon\colon\Decr(X)\to X$ defined using the discarded face maps, so that $\epsilon_n\colon\Decr(X)_n\to X_n$ is $d_{n+1}$; and a canonical map $\delta\colon\Decr(X)\to\Decr(\Decr(X))$ defined using the discarded degeneracy maps, so that $\delta_n\colon \Decr(X)_n\to\Decr(\Decr(X))_n$ is $s_{n+1}$. These maps $\delta$ and $\epsilon$ define the comultiplication and counit of the comonad.

Similarly, the left decalage $\Decl(X)$ of $X$ is the augmented simplicial set
$$\xymatrix{
\llap{\ldots~~}  X_3 \ar@<3ex>[r]|{d_1} \ar[r]|{d_2} \ar@<-3ex>[r]|{d_3} &
X_2 \ar@<1.5ex>[l]|{s_2} \ar@<-1.5ex>[l]|{s_1}  \ar@<1.5ex>[r]|{d_1} \ar@<-1.5ex>[r]|{d_2} &  X_1 \ar[l]|{s_1} \ar[r]^-{d_1} & X_0}$$
obtained by discarding $X_{-1}$ and the first face and degeneracy map in each degree.

We have described the decalage comonads for simplicial and augmented simplicial sets, but in much the same way there are decalage comonads $\Decr$ and $\Decl$ on the categories $[\Delta\op, P]$ and $[\bbm\op,P]$ of simplicial and augmented simplicial objects in $P$ for any category $P$, although in general there will no longer be a monoidal structure with respect to which decalage is given by an internal hom.






\subsection{Duplicial structure} \label{sect:duplicial}

Here we recall the definition of duplicial structure, and give a reformulation using the decalage comonads. As stated already in the introduction, a duplicial object in a category is a simplicial object $X$, equipped with a map $t_n\colon X_n\to X_n$ for each $n \geqslant 0$, subject to various conditions which we now state explicitly.
\begin{align}
\  d_i t_{n+1} &=
                \begin{cases}
                  t_n d_{i-1} & \text{if $1\le i\le n+1$;} \\
                  d_{n+1} & \text{if $i=0$;}
                \end{cases} \label{eq:1}\\
s_i t_n &=
          \begin{cases}
            t_{n+1} s_{i-1} & \text{if $1\le i\le n$;} \\
            t^2_{n+1} s_n & \text{if $i=0$.}
          \end{cases}\label{eq:2}
  \end{align}
There is also a formulation of this structure which uses an ``extra degeneracy map'' $s_{-1}\colon X_n\to X_{n+1}$ in each degree instead of the $t_n$; this $s_{-1}$ may be constructed as the composite $t_{n+1}s_n$. As in the introduction, $X$ is called \emph{paracyclic} if each $t_n$ is invertible, and \emph{cyclic} if additionally $t^{n+1}_n=1$.

The indexing category for cyclic structure is Connes' cyclic category $\Lambda$, which is a sort of wreath product of $\Delta$ and the cyclic groups. This is explained for example in \cite[Chapter~6]{Loday-cyclic}, where the more general notion of crossed simplicial group can also be found. This involves replacing the cyclic groups by some other family of groups indexed by the natural numbers, and equipped with suitable actions of $\Delta$ which allow the formation of the wreath product. The  indexing category for paracyclic structure can be obtained in this way on taking all the groups to be $\bbz$ \cite[Proposition~6.3.4 (c)]{Loday-cyclic}. Using the presentation for duplicial structure given above, it is straightforward to modify this argument to see that the indexing category $\bfk$ for duplicial structure is once again a wreath product, but this time by a ``crossed simplicial monoid'', involving the monoid $\bbn$ in each degree.

\begin{proposition}\label{prop:duplicial-decalage}
To give duplicial structure to a simplicial object $X$ is equivalently to give a simplicial map $t\colon \Decr X\to \Decl X$ making the following diagrams commute
$$\xymatrix{
\Decr X \ar[r]^t \ar[dr]_{\epsilon} & \Decl X \ar[d]^{\epsilon} \\ & X } \qquad
\xymatrix@C+0.3em{
\Decr X \ar[r]^t \ar[d]_{\delta} & \Decl X \ar[r]^{\delta} & \Decl^2X \\
\Decr^2 X \ar[r]_-{\Decr t} & \Decr\Decl X \ar@{=}[r] & \Decl\Decr X . \ar[u]_{\Decl t} }$$
\end{proposition}

\proof
The data of a simplicial map $t \colon \Decr X \rightarrow \Decl X$ comprises maps $t_n \colon X_n \rightarrow X_n$ for each $n > 0$ satisfying conditions. Compatibility of $t$ with face maps gives the cases where $i > 0$ of~\eqref{eq:1}, while those where $i=0$ are the compatibility condition with $\epsilon$. Likewise, compatibility of $t$ with degeneracy maps yields the cases $i,n > 0$ of~\eqref{eq:2}, while the cases where $n>0$ but $i=0$ are the compatibility condition with $\delta$.

The one thing which remains is to see that a map $t_0\colon X_0\to X_0$ satisfying~\eqref{eq:2} for $n = 0$ can be uniquely recovered from the remaining data and axioms. In order to have $s_0t_0 = t_1^2s_0$, we must have that $t_0 = d_0s_0t_0 = d_0t_1^2s_0 = d_1t_1s_0$. So we just need to check that, defining $t_0$ in this way, it satisfies the required relations; but this is indeed the case as the following calculations show:
\begin{gather*}
 (d_1 t_1 s_0) d_0 = d_1 t_1 d_0 s_1 = d_1 d_1 t_2 s_1 = d_1 d_2 t_2 s_1 = d_1 t_1 d_1 s_1 = d_1 t_1\text{ ; and} \\
s_0  (d_1 t_1 s_0) = d_2 s_0 t_1 s_0 = d_2 t^2_2 s_1 s_0 = t_1 d_1 t_2 s_1 s_0 = t^2_1 d_0 s_1 s_0 = t^2_1 d_0 s_0 s_0 = t^2_1 s_0\,.
\end{gather*}
\endproof

\subsection{The B{\"ohm}-{\c S}tefan construction}\label{sec:bohm-c-stefan}
We now describe the construction in~\cite{BohmStefan1,BohmStefan2}. The original formulation involves monads and coduplicial structure, but we work dually with comonads so as to obtain duplicial structure. Let $A$ and $P$ be categories, and suppose that we have a comonad $(g,\delta, \epsilon)$ on $A$ and a functor $f \colon A \to P$. As explained in Section~\ref{sect:bar}, we obtain from $g$ a functor $A \rightarrow [\bbm\op, A]$ sending each object to its bar resolution with respect to $g$, and post-composing with $f$ yields a functor $f^g \colon A \to [\bbm\op, P]$. Explicitly, $f^g$ takes $x$ in $A$ to the augmented simplicial object $f^g(x)$ with
$f^g(x)_n = fg^{n+1}x$ and with face and degeneracy maps:
\begin{align*}
 	d_i &=fg^i \epsilon 
 	g^{n-i} x \colon f^g(x)_n \rightarrow f^g(x)_{n-1} \\ \text{and} \qquad 
        s_j &= fg^j \delta g^{n-j} x \colon f^g(x)_n \rightarrow f^g(x)_{n+1}\rlap{ .}
\end{align*}
The basic construction of~\cite{BohmStefan1} uses additional data to equip objects of the form $f^g(x)$ with duplicial structure. We suppose given another comonad $h$ on $A$, and a \emph{distributive law}~\cite{Beck:dist-law} $\lambda \colon gh \to hg$---a natural transformation satisfying four axioms relating it to the comonad structures. We suppose moreover that the functor $f \colon A \rightarrow P$ is equipped with a natural transformation $\varphi \colon fh \rightarrow fg$ rendering commutative the diagrams:
\begin{equation}
  \label{eq:6}
  \vcenter{\hbox{\xymatrix @R1pc @C1pc@C+0.5em {
    fh \ar[rrr]^{\varphi} \ar[d]_{f\delta} &&& fg \ar[d]^{f\delta} && fh \ar[rr]^{\varphi} \ar[dr]_{f\epsilon} && fg\rlap{ .} \ar[dl]^{f\epsilon} \\
    fh^2 \ar[r]_{\varphi h} & fgh \ar[r]_{f\lambda} & fhg \ar[r]_{\varphi g} & fg^2 && & f
  }}}
\end{equation}
This was called \emph{left $\lambda$-coalgebra} structure on $f$ in \cite{Slevin2}, and the totality $(A,P,g,h,f,\lambda,\varphi)$
of the structure considered so far was called an {\em admissible septuple} in \cite{BohmStefan1}. Finally, we assume given an object $x \in A$ equipped with a map $\xi\colon gx\to hx$ rendering commutative:
\begin{equation}
  \label{eq:5}
  \vcenter{\hbox{\xymatrix @R1pc @C1pc@C+0.5em {
    gx \ar[rrr]^{\xi} \ar[d]_{\delta x} &&& hx \ar[d]^{\delta x} && gx \ar[rr]^{\xi} \ar[dr]_{\epsilon x} && hx \ar[dl]^{\epsilon x} \\
    g^2x \ar[r]_{g\xi} & ghx \ar[r]_{\lambda x} & hgx \ar[r]_{h\xi} & h^2x && & x
  }}}
\end{equation}
This was called {\em right $\lambda$-coalgebra} structure in \cite{Slevin1}, and a ``transposition map'' in~\cite{BohmStefan1}, though the notion itself goes back to~\cite{Burroni-NDA}. Under these assumptions, it was shown in~\cite{BohmStefan1} that the simplicial object $f^g(x)$ admits a duplicial structure. The duplicial operator $t_n \colon f^g(x)_n \to f^g(x)_n$ is given by the composite
$$
\xymatrix@=3em{
fg^{n+1}x \ar[r]^-{fg^n\xi  x} & fg^nhx\ar[r]^-{f\lambda^n  x} & fhg^n x \ar[r]^-{\phi  {g^n x}} & fg^{n+1} x
}
$$
where the natural transformation $\lambda^n \colon g^n h \to hg^n$ denotes the composite
$$
\xymatrix@=1.5em{
g^n h \ar[rr]^-{g^{n-1} \lambda} && g^{n-1} h g \ar[rr]^-{g^{n-2} \lambda g} && g^{n-2} h g^2 \ar[r] & \cdots \ar[r] & ghg^{n-1} \ar[rr]^-{\lambda g^{n-1}} && hg^n.}
$$
In~\cite{BohmStefan1}, this construction was used to obtain, among other things,  the cyclic cohomology and homology of bialgebroids. 

There is an automorphism $\Phi \colon [\bbm\op, P] \to [\bbm\op, P]$ induced by the automorphism in Remark~\ref{rmk:Delta-rev}, that maps a simplicial object $X$ to the simplicial object associated to $X$, obtained by reversing the order of all face and degeneracy maps. In~\cite{Slevin2} it is explained that $\Phi f^h ({x})$ is duplicial, and that there are two duplicial maps
$$
\xymatrix{
f^g({x}) \ar[r]^-R & \Phi f^h({x}),} \qquad \xymatrix{\Phi f^h({x}) \ar[r]^-L & f^g({x}),
}
$$
defined by iteration of $\phi$, respectively $\xi$, which are mutual inverses if and only if both objects are cyclic.

\begin{remark} 
The distributive laws $\lambda\colon gh\to hg$ considered above are the objects of the 2-category $\Mnd_*^*(\Mnd_*(\Cat))$  of \cite{ftm}. It was observed in \cite{Slevin1, Slevin2} that the morphisms of this 2-category act on right $\lambda$-coalgebras $(x,\xi)$ and on left $\lambda$-coalgebras $(f,\phi)$, thereby giving rise to more examples of duplicial objects; in the lax Hochschild theory developed below, this observation will become the functoriality of homology and cohomology. One of the starting points for this paper was the observation that $f^g$ can be seen as a morphism of this 2-category from $(A,g,h)$ to $([\bbm\op,P],\Decr,\Decl)$, with the action of this morphism on the left $\lambda$-coalgebra $(x,\xi)$ yielding $(f^g(x),t)$, the duplicial object of B\"ohm-\c Stefan. 
\end{remark}


\subsection{Zeroth Hochschild homology and cohomology}\label{sect:Hochschild}

Let $A$ be a ring, and $X$ a bimodule over $A$. There is an induced simplicial abelian group part of which looks like
$$\xymatrix{ \cdots &
A\ox A\ox X \ar@<1.5ex>[r]|-{d_0} \ar[r]|-{d_1} \ar@<-1.5ex>[r]|-{d_2} & A\ox X \ar@<1.5ex>[r]|-{d_0} \ar@<-1.5ex>[r]|-{d_1} & X \ar[l]|-{s_0} }$$
with the maps given as follows:
\begin{align*}
  d_0(a\ox x) &= xa & d_0(a\ox b\ox x) &= b\ox xa \\
  d_1(a\ox x) &= ax & d_1(a\ox b\ox x) &= ab\ox x \\
  s_0(x) &= 1\ox x & d_2(a\ox b\ox x) &= a\ox bx\rlap{ ,}
\end{align*}
and which is defined analogously in higher degrees.  
We call this simplicial object the \emph{Hochschild complex} of $X$, although often that name refers to the corresponding (normalized or otherwise) chain complex.

The \emph{zeroth homology of $A$ with coefficients in $X$} is the colimit $H_0(A,X)$ of this diagram, which can more simply be computed as the coequalizer of the two maps $A\ox X\rightrightarrows X$; more explicitly still, this is the quotient of $X$ by the subgroup generated by all elements of the form $ax-xa$

Dually there is a cosimplicial object part of which looks like
$$\xymatrix{
X \ar@<1.5ex>[r]|-{\delta_0} \ar@<-1.5ex>[r]|-{\delta_1} & [A,X] \ar[l]|-{\sigma_0} \ar@<1.5ex>[r]|-{\delta_0} \ar[r]|-{\delta_1} \ar@<-1.5ex>[r]|-{\delta_2} & [A \otimes A,X] & \cdots}$$
with the maps given as follows
\begin{align*}
  \delta_0(x)(a) &= xa & \delta_0(f)(a\ox b) &= f(a)b \\
\delta_1(x)(a) &= ax & \delta_1(f)(a\ox b) &= f(ab) \\
\sigma_0(f) &= f(1) & \delta_2(f)(a\ox b) &= a f(b)\rlap{ ,}
\end{align*}
and now the \emph{zeroth Hochschild cohomology of $A$ with coefficients in $X$} is the limit $H^0(A,X)$ (really an equalizer) of this diagram, given explicitly by the subgroup of $X$ consisting of those $x$ for which $xa=ax$ for all $a\in A$.

\subsection{Universality of zeroth Hochschild homology and cohomology}
\label{sec:univ-zeroth-hochsch}
Both $H^0(A,X)$ and $H_0(A,X)$ have universal characterisations. For any $A$-bimodule $X$ and any abelian group $P$, there is an induced bimodule structure on $[X,P]$ given by $(af)(x) = f(xa)$ and $(fa)(x) = f(ax)$, and this construction gives a functor $[X, \mathord -] \colon \mathbf{Ab} \rightarrow A\bMod A$. In particular we may take $X = A$ with its regular left and right actions.

\begin{proposition}\label{prop:4}
The functor $[A, \mathord -] \colon \mathbf{Ab} \rightarrow A \bMod A$ has a left adjoint sending an $A$-bimodule $X$ to $H_0(A,X)$.
\end{proposition}
Similarly, there is for any $A$-bimodule $X$ and abelian group $P$ an induced bimodule structure on $X \otimes P$ given by $a(x \otimes p) = ax \otimes p$ and $(x \otimes p)a = xa \otimes p$, and this gives a functor $X \otimes (\mathord -) \colon \mathbf{Ab} \rightarrow A \bMod A$. Considering again the case $X = A$, we have:
\begin{proposition}
\label{prop:2}
The functor $A \otimes (\mathord -) \colon \mathbf{Ab} \rightarrow A \bMod A$ has a right adjoint sending an $A$-bimodule $X$ to $H^0(A,X)$.
\end{proposition}

\section{Bimodules} \label{sect:bimodules}

We described above the Hochschild complex of a ring $A$ with coefficients in an $A$-bimodule. A ring is the same thing as a monoid in the monoidal category $\Ab$ of abelian groups, and more generally the Hochschild complex and the zeroth homology and cohomology can be constructed if $A$ is a monoid in a suitable symmetric monoidal closed category \cv. In particular, we could do this for the cartesian closed category \Cat. But \Cat is in fact a 2-category, which opens the way to consider \emph{lax} variants of the theory, and it is such a variant that we shall now present. While it would be possible to develop this theory in the context of a general symmetric monoidal closed bicategory \cv, it is only the case $\cv = \Cat$ which will need, and so we restrict ourselves to that. The first step, carried out in this section, is to describe in detail the notion of bimodule that will play the role of coefficient object for our lax homology and cohomology.

\subsection{Monoids}

A monoid in \Cat is precisely a strict monoidal category. It is not particularly difficult to adapt the theory that follows to deal with non-strict monoidal categories, but we do not need this extra generality, and feel that the complications that it causes might distract from the story we wish to tell. It is probably also possible to extend the theory to deal with skew monoidal categories \cite{Szlachanyi-skew,skew}, although we have not checked this in detail.

We shall therefore consider a strict monoidal category $(A,m,i)$. We shall write $a\ox b$ or sometimes just $ab$ for the image under the tensor functor $m\colon A\x A\to A$ of a pair $(a,b)$.

\subsection{Modules}

Next we need a notion of module over $A$. There is a well-developed (pseudo) notion of an action of a monoidal category on a category, sometimes called an actegory. Here, however, we deal only with the strict case, which does not use the 2-category structure of \Cat; once again it would not be difficult to extend our theory to deal with pseudo (or possibly skew) actions, but this is not needed for our applications so we have not done so. To give a strict left action of $A$ on a category $X$ is equivalently to give a strict monoidal functor from $A$ to the strict monoidal category $\End(X)$ of endofunctors of $X$. The image under the corresponding functor $\alpha\colon A\x X\to X$ of an object $(a,x)$ will be written $ax$. Similarly there are (strict) right actions involving functors $\beta\colon X\x A\to X\colon (x,a)\mapsto xa$ satisfying strict associativity and unit conditions.

In fact we shall also make use of a slightly more general notion. It is possible to consider actions of monoids not just on sets, but also on objects of other categories; in the same way it is possible to consider actions of monoidal categories on objects of other 2-categories. If $X$ is an object of a 2-category \ck, then an action of $A$ on $X$ will be a strict monoidal functor from $A$ to the strict monoidal category $\ck(X,X)$ of endomorphisms of $X$.

If the 2-category \ck admits copowers, then there is an equivalent formulation as follows. Recall that the copower of an object $X$ by a category $P$ is an object $P\cdot X$ equipped with isomorphisms of categories
$$\ck(P\cdot X,Y) \cong \Cat(P,\ck(X,Y))$$
2-natural in the variable $Y\in\ck$. If \ck has all copowers, then there are 2-natural isomorphisms $(P\x Q)\cdot X\cong P\cdot(Q\cdot X)$ and $1\cdot X\cong X$. In this case, a strict (left) action of $A$ on $X$ is equivalently a morphism $\alpha\colon A\cdot X\to X$ in $\ck$ for which the diagrams
$$\xymatrix{
(A\x A)\cdot X \ar[rr]^{m\cdot 1} \ar[d] && A\cdot X \ar[d]^{\alpha} & 1\cdot X \ar[r]^{i\cdot 1} \ar[dr] & A\cdot X \ar[d]^{\alpha} \\
A\cdot(A\cdot X) \ar[r]_{1\cdot \alpha} & A\cdot X \ar[r]_{\alpha} & X & & X }$$
commute, where the un-named maps are the isomorphisms just described. (There are also still more general notions of action of $A$: see \cite[Section 2]{property}.)

Note that the 2-category \Cat admits copowers, with $A\cdot X$ given by the cartesian product $A\x X$, so that in this case our more general notion of action of $A$ on $X \in \Cat$ reduces to the initial one.

\begin{example}
Our running example throughout this section and the next will take $A$ to be the strict monoidal category $\bbm\op$; it is this example which will be used to explain the B\"ohm-\c Stefan construction. Since a strict monoidal functor $\bbm\op\to\ck(X,X)$ is precisely a comonoid in $\ck(X,X)$, a left $\bbm\op$-module is a comonad in the 2-category \ck, in the sense of \cite{ftm}. On the other hand, a right $\bbm\op$-module is also just  a comonad in \ck, as follows from Remark~\ref{rmk:Delta-rev}.

In the case $\ck=\Cat$, a comonad in $\Cat$ is a category $X$ equipped with a comonad $g$. For an object $n$ of $\bbm\op$ and an object $x\in X$, the value $nx$ of the corresponding left $\bbm\op$-action is given by $g^nx$.
\end{example}

\subsection{Morphisms of modules}

When it comes to morphisms of modules, once again there is a question of how lax they should be, and this time we deviate from the completely strict situation. If $X$ and $Y$ are (strict, as ever) left $A$-modules in \Cat, we define a {\em lax $A$-morphism} to be a functor $p\colon X\to Y$, equipped with a natural transformation
$$\xymatrix{
A\x X \ar[r]^{1\x p}_{~}="1" \ar[d]_{\alpha} & A\x Y \ar[d]^{\alpha} \\
X \ar[r]_{p}^{~}="2" & Y
\ar@{=>}"1";"2"^{\rho}
}$$
whose components have the form
$$\xymatrix{
a.p(x) \ar[r]^{\rho_{a,x}} & p(ax) }$$
for  $a\in A$ and $x\in X$, and which satisfy two coherence conditions. The first asks that $\rho_{i,x}\colon p(x)=i.p(x)\to p(ix)=p(x)$ is the identity. The second ask that the composite
$$\xymatrix{
ab.p(x) \ar[r]^{a\rho_{b,x}} & a.p(bx) \ar[r]^{\rho_{a,bx}} & p(abx) }$$
be equal to $\rho_{ab,x}$. Often we omit the subscripts and simply write $\rho$ for $\rho_{a,x}$. When $\rho$ is an  identity, we say that the $A$-morphism is {\em strict}.

For actions on objects of a general 2-category \ck given by strict monoidal functors $A\to\ck(X,X)$ and $A\to \ck(Y,Y)$, a lax $A$-morphism will be a morphism $p\colon X\to Y$ in \ck together with a natural transformation
$$\xymatrix{
A \ar[r]_{~}="1" \ar[d] & \ck(X,X) \ar[d]^{\ck(X,p)} \\ \ck(Y,Y) \ar[r]_{\ck(p,Y)}^{~}="2" & \ck(X,Y)
\ar@{=>}"2";"1" }$$
satisfying an associativity and a unit axiom generalizing those above. If \ck admits copowers, then the natural transformation displayed above determines and is determined by a 2-cell
$$\xymatrix{
A\cdot X \ar[r]^{1\cdot p}_{~}="1" \ar[d]_{\alpha} & A\cdot Y \ar[d]^{\alpha} \\
X \ar[r]_{p}^{~}="2" & Y
\ar@{=>}"1";"2"^{\rho}
}$$
in \ck, satisfying associativity and unit conditions.

If $(p,\rho)$ and $(p',\rho')$ are lax $A$-morphisms from $X$ to $Y$, an $A$-transformation from $(p,\rho)$ to $(p',\rho')$ is a 2-cell $\tau\colon p\to p'$ satisfying the evident compatibility condition; in the case $\ck=\Cat$, this says that the diagram
$$\xymatrix{
a.px \ar[d]_{\rho_{a,x}} \ar[r]^{1.\tau} & a.p'x \ar[d]^{\rho'_{a,x}} \\
p(ax) \ar[r]_{\tau} & p'(ax) }$$
commutes for all objects $a\in A$ and $x\in X$.

There is a 2-category $A\lMod$ whose objects are the $A$-modules (in \Cat), whose morphisms are the lax $A$-morphisms, and whose 2-cells are the $A$-transformations. This 2-category admits copowers, with $B\cdot X$ given by the category $X\x B$ equipped with the action $\alpha\x 1\colon A\x X\x B\to X\x B$, where $\alpha\colon A\x X\to X$ is the action on $X$.

\begin{example}
In the case $A=\bbm\op$, we saw that an $A$-module was precisely a category $X$ equipped with a comonad $g$. A lax $A$-morphism is what was called a {\em comonad opfunctor} in \cite{ftm}, and indeed $\bbm\op\lMod$ is the 2-category called $\Mnd^*_*(\Cat^*_*)$ in that paper.
\end{example}

\subsection{Bimodules}
As usual, a bimodule is an object which is both a left and right module with suitable compatibility between the two actions. Although our notion of action is strict, the compatibility between the actions will not be. A succinct definition of $A$-bimodule is: an object of $A\lMod$ equipped with a right $A$-module structure, but we can also spell out what this means.

First of all, there is a category $X$ with a strict left action $\alpha\colon A\x X\to X$.
The right action involves a functor $\beta\colon X\x A\to X$ defining a strict right action, but this should be not just a functor, but a lax $A$-module morphism $A\cdot X\to X$. This lax $A$-morphism structure consists of maps
$$\xymatrix{
a(xb) \ar[r]^{\lambda_{a,x,b}} & (ax)b }$$
natural in the variables $a\in A, x\in X, b\in A$, and making each diagram
$$\xymatrix{
aa'(xb) \ar[r]^{1\lambda_{a',x,b}} \ar[dr]_{\lambda_{aa',x,b}} & a((a'x)b) \ar[d]^{\lambda_{a,a'x,b}} & i(xb) \ar[dr]_{\lambda_{i,x,b}} \ar@{=}[r] & xb \ar@{=}[d] \\
& ((aa')x)b & & (ix)b }$$
commute. Finally, the associative and unit laws required for the right action $\beta\colon X\x A\to X$ should hold not just as equations between functors, but as equations between lax $A$-morphisms. Explicitly, this means that each diagram
$$\xymatrix{
a(xbb') \ar[r]^{\lambda_{a,xb,b'}} \ar[dr]_{\lambda_{a,x,bb'}} & (a(xb))b' \ar[d]^{\lambda_{a,x,b}1} & a(xi) \ar@{=}[r] \ar[dr]_{\lambda_{a,x,i}} & ax \ar@{=}[d] \\
& (ax)(bb') && (ax)i }$$
should commute.

\begin{example}
Returning to our running example $A=\bbm\op$, we have already seen that the 2-category $A\lMod$ is just Street's 2-category $\Mnd^*_*(\Cat^*_*)$ of comonads and comonad opfunctors, and that a right $\bbm\op$-action in a 2-category is a comonad in that 2-category. So an $A$-bimodule will be a comonad in $\Mnd^*_*(\Cat^*_*)$, which as explained in \cite{ftm} amounts to a category $X$ equipped with comonads $g$ and $h$ and a distributive law $\lambda\colon gh\to hg$ between them.
\end{example}

\subsection{Morphisms of bimodules}

While our morphisms of left modules are lax, we shall consider only strict morphisms of right modules, but once again these should be defined relative to the 2-category $A\lMod$.
This means that a morphism $(X,\alpha,\beta)\to(Y,\alpha,\beta)$ of bimodules will be a lax $A$-morphism $(p,\rho)\colon (X,\alpha)\to(Y,\alpha)$ of the underlying left modules, for which the diagram
$$\xymatrix{
X\x A \ar[r]^{p\x 1} \ar[d]_\beta & Y\x A \ar[d]^\beta \\ X \ar[r]_p & Y }$$
of categories and functors commutes, and for which moreover the diagram
\begin{equation}\label{eq:bimodule-map}
\xymatrix @R1pc@C+0.3em {
& a.(px.b) \ar[r]^{\lambda_{a,px,b}} & (a.px).b \ar[dr]^{\rho_{a,x}.1} \\
a.p(xb) \ar@{=}[ur] \ar[dr]_{\rho_{a,xb}} & & & p(ax).b \\
& p(a(xb)) \ar[r]_{p(\lambda_{a,x,b})} & p((ax)b) \ar@{=}[ur] }
\end{equation}
commutes for all $a,b\in A$ and $x\in X$.

The bimodules and their morphisms constitute the objects and morphisms of a 2-category $A\bMod A$; a 2-cell $(p,\rho)\to (p',\rho')$ is a natural transformation $\tau\colon p\to p'$ which is a 2-cell relative to both the left and right actions.

\begin{example}\label{ex:[A,-]}
 For an $A$-bimodule $X$ and an arbitrary category $P$, the functor category $[X,P]$ has left and right actions of $A$, given by $(af)(x)=f(xa)$ and $(fa)(x)=f(ax)$, and these define a bimodule structure on $[X,P]$. This forms the object part of a 2-functor $[X,-]\colon \Cat\to A\bMod A$. We shall be particularly interested in the case where $X$ is $A$ with its standard bimodule structure; in this case, since the left and right actions on $A$ are strictly compatible, so too are those on $[A,P]$.
\end{example}

\begin{example}\label{ex:-xA}
  Dually, for an $A$-bimodule $X$ and an arbitrary category $P$, the product category $P\x X$ has left and right actions inherited from $X$, and this forms the object part of a 2-functor $(\mathord -)\x X\colon \Cat\to A\bMod A$.
\end{example}

\section{Lax cohomology and homology}\label{laxhom}

\subsection{The Hochschild complex}

Let $A$ be a strict monoidal category and $X$ a bimodule over $A$, in the sense of the previous section. Then we can define maps  
\begin{equation}
  \label{eq:4}
\xymatrix{ \cdots &
A\x A\x X \ar@<1.5ex>[r]|-{d_0} \ar[r]|-{d_1} \ar@<-1.5ex>[r]|-{d_2} & A\x X \ar@<1.5ex>[r]|-{d_0} \ar@<-1.5ex>[r]|-{d_1} & X \ar[l]|-{s_0} }
\end{equation}
exactly as in Section~\ref{sect:Hochschild}, except that, because of the lax compatibility between the actions, the simplicial identity $d_1d_0=d_0d_2$ no longer holds; instead, there is a natural transformation $\lambda\colon d_1d_0\to d_0d_2$ whose component at an object $(b,a,x)\in A\x A\x X$ is the map $\lambda_{a,x,b}\colon a(xb)\to (ax)b$. If we included more of the maps from the Hochschild complex of Section~\ref{sect:Hochschild}, then further equalities would be replaced by natural transformations, and using these we could express the various coherence conditions on $\lambda$ which appear in the definition of $A$-bimodule.

Similarly, there are maps
\begin{equation}
  \label{eq:3}
  \xymatrix{
  X \ar@<1.5ex>[r]|-{\delta_0} \ar@<-1.5ex>[r]|-{\delta_1} & [A,X]
  \ar[l]|-{\sigma_0} \ar@<1.5ex>[r]|-{\delta_0} \ar[r]|-{\delta_1}
  \ar@<-1.5ex>[r]|-{\delta_2} & [A \x A,X] & \cdots}
\end{equation}
defined as in Section~\ref{sect:Hochschild} once again; this time the
cosimplicial identity $\delta_2\delta_0 = \delta_1 \delta_0$ becomes a
natural transformation $\delta_2\delta_0\to \delta_1\delta_0$, whose
components are once again induced by the lax compatibilities
$\lambda_{a,x,b}$.

\subsection{Cohomology}

In Section~\ref{sect:Hochschild}, the zeroth Hochschild cohomology group $H^0(A,X)$ of a bimodule over a ring was defined as the equalizer of the maps $\delta_0,\delta_1\colon X \rightrightarrows [A,X]$. In the case of the lax cohomology of a bimodule over a strict monoidal category $A$, we define the zeroth Hochschild cohomology $H^0(A,X)$ by taking a ``lax version'' of an equalizer, involving all of the data displayed in~\eqref{eq:3}, called a \emph{lax descent object}; this is a mild variant~\cite{codescent} of a notion introduced in~\cite{FibBicCorr}. Interpreting this for~\eqref{eq:3} yields that $H^0(A,X)$ is the universal category $Y$ equipped with a functor $y\colon Y\to X$ and a natural transformation $\xi\colon \delta_1 y\to \delta_0 y$ such that $\sigma_0\xi\colon x=\sigma_0\delta_1y\to \sigma_0\delta_0 y=y$ is the identity and the diagram
$$\xymatrix @R1pc @C1pc {
& \delta_2\delta_0y \ar[r]^{\lambda y} & \delta_0\delta_1y \ar[dr]^{\delta_0\xi} \\
\delta_2\delta_1y \ar[ur]^{\delta_2\xi} \ar@{=}[dr] &&& \delta_0\delta_0y \\
& \delta_1\delta_1y \ar[r]_{\delta_1\xi} & \delta_1\delta_0y \ar@{=}[ur] }$$
commutes. Explicitly, an object of $H^0(A,X)$ is an object $x\in X$ equipped with maps $\xi_a\colon ax\to xa$ natural in $a\in A$, and satisfying $\xi_i=1$ as well as the cocycle condition asserting that the diagram
$$\xymatrix @R1pc @C1pc  {
& a(xb) \ar[r]^{\lambda x} & (ax)b \ar[dr]^{\xi b} \\
a(bx) \ar[ur]^{a\xi} \ar[dr] &&& (xa)b \\
& (ab)x \ar[r]_{\xi} & x(ab) \ar[ur] }$$
commutes for all $a,b \in A$.

\begin{example}
In the case of classical Hochschild cohomology, for a ring $A$ the zeroth cohomology group $H^0(A,A)$ is the centre of the ring; similarly for a strict monoidal category $A$, the lax cohomology $H^0(A,A)$ is the {\em lax centre} of $A$ in the sense of \cite{LaxCentre}, originally introduced in \cite{Schauenburg-DualsAndDoubles} with the name {\em weak centre}.
\end{example}

\begin{example}\label{ex:0}
  Consider our running example of $A=\bbm\op$, so that an $A$-bimodule $X$ is a category equipped with comonads $g$ and $h$ and a distributive law $\lambda\colon gh\to hg$. Explicit calculation shows that an object of $H^0(A,X)$ is an object $x\in X$ equipped with a map $\xi\colon gx\to hx$ making the diagrams~\eqref{eq:5} commute: so we re-find the notion of right $\lambda$-coalgebra of Section~\ref{sec:bohm-c-stefan}.
\end{example}

The next result justifies the definition of the lax cohomology
$H^0(A,X)$ analogously to Proposition~\ref{prop:2} for the usual Hochschild cohomology.
\begin{theorem}
The 2-functor $(\mathord -)\x A\colon \Cat\to A\bMod A$ has a right adjoint sending an $A$-bimodule $X$ to $H^0(A,X)$.
\end{theorem}

\proof
Let $X$ be an $A$-bimodule and $P$ a category. To give a (strict) right $A$-module morphism $p\colon P\x A\to X$ is equivalently to give a functor $f\colon P\to X$; here $f(y)=p(y,1)$ and $p(y,a)=f(y)a$. To enrich such a morphism of modules into a morphism $(p,\rho)$ of bimodules, we should give suitably natural and coherent maps $\rho_{a,y,b}\colon a.p(y,b)\to p(y,ab)$ for all $a\in A$ and $(y,b)\in P\x A$. By the compatibility condition \eqref{eq:bimodule-map}, the map $\rho_{a,y,b}$ can be constructed as
$$\xymatrix @C3pc {
{}\llap{$ap(y,b)= $ } a(p(y,1)b) \ar[r]^{\lambda_{a,p(y,1),b}} & (ap(y,1))b \ar[r]^{\rho_{a,y,1}1} & p(y,a)b \rlap{ $=p(y,ab)$}
}$$
and so the general $\rho$ will be determined by those of the form $\rho_{a,y,1}$, and these have the form $\xi_{a,y}\colon af(y)\to f(y)a$. The unit condition asserting that each $\rho_{1,y,b}$ is the identity says that $\xi_{1,y}$ is the identity. The cocycle condition on the $\rho$ is equivalent to the cocycle condition asserting that $\xi_{a,y}$ makes each $f(y)$ into an object of $H^0(A,X)$. Naturality of $\xi_{a,y}$ in $y$ implies that for each morphism $\psi\colon y\to y'$ in $P$, the map $f(\psi)$ defines a morphism $(f(y),\xi_{a,y})\to (f(y'),\xi_{a,y'})$ in $H^0(A,X)$.

This gives the desired bijection between bimodule morphisms $P\x A\to X$ and functors $P\to H^0(A,X)$; it is straightforward to check that this carries over to 2-cells, and so defines an isomorphism of categories
$$A\bMod A(P\x A,X)\cong \Cat(P,H^0(A,X))$$
exhibiting $H^0(A,X)$ as the value at $X$ of a right adjoint to $(\mathord -) \x A$.
\endproof

\subsection{Homology}

In Section~\ref{sect:Hochschild} the zeroth Hochschild homology group was defined as the coequalizer of the maps $d_0,d_1\colon A\ox X\rightrightarrows X$. For lax homology, we define $H_0(A,X)$ of an $A$-bimodule $X$ to be the \emph{lax codescent object} of the data displayed in~\eqref{eq:4}. Lax codescent objects are the colimit notion corresponding to the lax descent objects used to define lax cohomology.

Spelling this out, $H_0(A,X)$ is the universal category $Y$ equipped with a functor $f\colon X\to Y$ and a natural transformation $\phi\colon fd_0\to fd_1$ satisfying the normalization condition $\phi s_0=1$ and the cocycle condition
$$\xymatrix @R1pc @C1pc @C+0.5em{
& fd_1d_0 \ar[r]^{f\lambda} & fd_0d_2 \ar[dr]^{\phi d_2} \\
fd_0d_0 \ar[ur]^{\phi d_0} \ar@{=}[dr] &&& fd_1d_2\rlap{ .} \\
& fd_0d_1 \ar[r]_{\phi d_1} & fd_1d_1 \ar@{=}[ur] }$$
Explicitly, $H_0(A,X)$ is obtained from $X$ by adjoining morphisms $xa\to ax$ satisfying naturality conditions in both variables, with $xi\to ix$ required to be the identity, and obeying the cocycle condition which requires the diagram
$$\xymatrix @R1pc @C1pc@C+1em {
& b(xa) \ar[r]^{\lambda} & (bx)a \ar[dr]^{\phi_{bx,a}} \\
(xa)b \ar[dr] \ar[ur]^{\phi_{xa,b}} &&& a(bx) \\
& x(ab) \ar[r]_{\phi_{x,ab}} & (ab)x \ar[ur] }$$
to commute.

\begin{example}\label{ex:1}
Let $A=\bbm\op$, and let $X$ have $A$-bimodule structure corresponding to comonads $g$ and $h$ and a distributive law $\lambda\colon gh\to hg$. By the defining universal property of the category $H_0(A,X)$, giving a functor $H_0(A,X)\to P$ is the same as giving a functor $f \colon A \rightarrow P$ and natural transformation $\varphi \colon fh \rightarrow fg$ making the diagrams~\eqref{eq:6} commute: so we re-find the notion of left $\lambda$-coalgebra from Section~\ref{sec:bohm-c-stefan}.
\end{example}
\begin{example}\label{ex:2}
Again with $A = \bbm\op$, the ``regular'' $A$-bimodule structure on $A$ corresponds to the two decalage comonads equipped with the identity distributive law between them. The full subcategory of $\bbm\op$ given by the non-empty finite ordinals is a sub-bimodule; since it is also isomorphic to $\Delta\op$, there is an induced bimodule structure on $\Delta\op$. By the preceding example and the description of duplicial structure given in Proposition~\ref{prop:duplicial-decalage}, a functor $H_0(\bbm\op,\Delta\op) \rightarrow P$ is precisely a duplicial object in $P$, so that $H_0(\bbm\op,\Delta\op)$ itself is the category $\mathbf{K}\op$ indexing duplicial structure. Similarly, a functor $H_0(\bbm\op,\bbm\op) \to P$ is an augmented duplicial object in $P$, and $H_0(\bbm\op,\bbm\op)$ is the category indexing augmented duplicial structure.
\end{example}

Just as before, the lax zeroth Hochschild homology has a universal characterisation paralleling Proposition~\ref{prop:4}.

\begin{theorem}\label{thm:homology-adjunction}
The 2-functor $[A,-]\colon \Cat\to A\bMod A$ has a left adjoint sending an $A$-bimodule $X$ to $H_0(A,X)$.
\end{theorem}

\proof
Let $X$ be an $A$-bimodule and $P$ a category. Just as in the classical case, to give a (strict) morphism of right $A$-modules $p\colon X\to [A,P]$ is equivalently to give a morphism $f\colon X\to P$, with $f(x)=p(x)(1)$ and $p(x)(a)=f(xa)$. In order to enrich such a $p$ into a morphism $(p,\rho)\colon X\to [A,P]$ of bimodules, we should give a suitably coherent map $\rho_{a,x}\colon a.p(x)\to p(ax)$ in $[A,P]$ for all $a\in A$ and $x\in X$. Thus for $b\in A$ we should give
$$\xymatrix{ {}\llap{$f(x(ba))=p(x)(ba)= $} (a.p(x))b\ar[r]^-{\rho_{a,x}(b)} &  p(ax)(b) \rlap{ $= f((ax)b).$ } }$$
Commutativity of \eqref{eq:bimodule-map} means that the general $\rho_{a,x}(b)$ is equal to the composite
$$\xymatrix{
a.p(xb) \ar[r]^{\rho_{a,xb}(1)} & p(a(xb)) \ar[r]^{p\lambda_{a,x,b}} & p((ax)b) \rlap{ $=p(ax)b$.} }$$
Thus $\rho$ is determined by the maps $\rho_{a,x}(1)\colon f(xa)\to f(ax)$, which we can regard as defining a natural transformation $\phi\colon fd_0\to fd_1$. The normalization condition asserting that $\rho_{1,x}$ is an identity now says that
$$\xymatrix{
{}\llap{$f=$ } fd_0s_0 \ar[r]^{\phi s_0} & fd_1s_0 \rlap{ $=f$} }$$
is an identity. The cocycle condition on the $\rho$ is equivalent to the cocycle condition on $\phi$, and so we have a bijection between bimodule morphisms $X\to [A,P]$ and functors $H_0(A,X)\to P$. It is straightforward to extend this to 2-cells, and so to obtain an isomorphism of categories
$$A\bMod A(X,[A,P]) \cong \Cat(H_0(A,X),P)$$
exhibiting $H_0(A,X)$ as the value at $X$ of a left adjoint to $[A,-]$.
\endproof

\subsection{The universal coefficients theorem and the cap product}

In this section we develop a few very simple ingredients of classical Hochschild theory in our lax context. The first of these is the universal coefficients theorem. In its more general forms this involves short exact sequences connecting homology and cohomology, but in degree zero it is particularly simple.

\begin{proposition}[Universal Coefficients Theorem] \label{prop:universal-coefficients}
 For any bimodule $X$ and category $P$ there is an isomorphism of categories
$$\Cat(H_0(A,X),P) \cong H^0(A,[X,P])$$
natural in $X$ and $P$.
\end{proposition}

\proof
By the universal property of $H_0(A,X)$ as a lax codescent object, an object of the left hand side amounts to a functor $f\colon X\to P$ equipped with a natural transformation $\phi\colon fd_0\to fd_1$ satisfying the normalization and cocycle conditions. But the functor $f$ can be seen as an object of $[X,P]$, while $\delta_0(f)\colon A\to [X,P]$ and $\delta_1(f)$ correspond under the adjunction $-\x A\dashv \Cat(A,-)$ to $fd_0\colon A\x X\to P$ and $fd_1$, so that to give $\phi\colon fd_0\to fd_1$ is equivalently to give $\xi\colon \delta_0(f)\to \delta_1(f)$. A straightforward calculation shows that the normalization and cocycle conditions for $\phi$ to make $f$ into a functor $H_0(A,X)\to P$ are equivalent to the normalization and cocycle conditions for $\xi$ to make $f$ into an object of $H^0(A,[X,P])$.

This proves that we have a bijection on objects; the case of morphisms is similar but easier, and is left to the reader.
\endproof

\begin{construction}[Cap product] \label{const:cap}
  For any $A$-bimodule $X$, we have the unit $\chi\colon X\to [A,H_0(A,X)]$ of the adjunction $H_0(A,-)\dashv [A,-]$ of Theorem~\ref{thm:homology-adjunction}. Applying the cohomology 2-functor $H^0(A,-)$ we obtain a functor
$$\xymatrix@C+1em{
H^0(A,X) \ar[r]^-{H^0(A,\chi)} & H^0(A,[A,H_0(A,X)]) }$$
and composing with the ``universal coefficients'' isomorphism $H^0(A,[A,P])\cong \Cat(H_0(A,A),P)$ of Proposition~\ref{prop:universal-coefficients}, we obtain a functor
$$\xymatrix{
H^0(A,X) \ar[r] & \Cat(H_0(A,A),H_0(A,X)) }$$
whose adjoint transpose
$$\xymatrix{
H^0(A,X) \x H_0(A,A) \ar[r] & H_0(A,X) }$$
can be seen as a special case of the cap product for our lax homology and cohomology. But we choose instead to transpose again to obtain a functor
$$\xymatrix{
H_0(A,A) \ar[r]^-{\BS} & \Cat(H^0(A,X),H_0(A,X)) }$$
which we call the {\em B\"ohm-\c Stefan} map.
\end{construction}

\begin{example}
  We now analyze this B\"ohm-\c Stefan map in the case of our running example. Suppose then that $A=\bbm\op$, and $X$ is an $A$-bimodule, with the bimodule structure corresponding to comonads $g$ and $h$ and a distributive law $\lambda\colon gh\to hg$. Let $p\colon H_0(A,X)\to P$ be an arbitrary functor, and let $y \in H^0(A,X)$.
As in Example~\ref{ex:1}, to give $p$ is equivalently to give a functor $f\colon X\to P$ equipped with left $\lambda$-coalgebra structure $\phi\colon fh\to fg$, while as in Example~\ref{ex:0}, to give $y$ is equally to give an object $x\in X$ equipped with right $\lambda$-coalgebra structure $\xi \colon gx\to hx$. There is now an induced functor
$$\xymatrix{
H_0(A,A) \ar[r]^-{\BS} & \Cat(H^0(A,X),H_0(A,X)) \ar[r]^-{\ev_y} & H_0(A,X) \ar[r]^-{ p} & P}$$
which by Example~\ref{ex:2} picks out an augmented duplicial object in $P$. This object is precisely the one constructed in~\cite{BohmStefan1} as recalled in Section~\ref{sec:bohm-c-stefan} above. This construction was generalized slightly in \cite{BohmStefan3}  to include right $\lambda$-coalgebra structures on arbitrary functors $Y \rightarrow X$, rather than just objects of $X$; in this case $y$ becomes a functor $Y\to H^0(A,X)$ and the composite
$$\xymatrix @C3pc {
H_0(A,A) \ar[r]^-{\BS} & \Cat(H^0(A,X), H_0(A,X)) \ar[r]^-{\Cat(y,p)} & \Cat(Y,P) }$$
defines an augmented  duplicial object in $\Cat(Y,P)$.
\end{example}

\section{Duplicial structure on nerves}

In this section we turn to our second main goal, which is to analyze duplicial structure on nerves of various sorts of categorical structures; specifically, on categories, on monoidal categories, and on bicategories.

A monoidal category can of course be seen as a one-object bicategory, and a category can be seen as a bicategory with no non-identity 2-cells, so in principle we could pass straight to the case of bicategories, and then merely read off the results for the other two cases, but instead we have chosen to do the case of categories first, as a sort of warm-up.

\subsection{Duplicial structure on categories}

The nerve functor from \Cat to \SSet is of course fully faithful, so that we may identify (small) categories with certain simplicial sets. It therefore makes sense to speak of duplicial structure borne by a category. The decalage comonads on \SSet restrict to \Cat, and so we may analyze duplicial structure on categories using Proposition~\ref{prop:duplicial-decalage}.

The right decalage comonad sends a category $C$ to the coproduct $\sum_x C/x$ over all objects $x\in C$ of the corresponding slice categories. The counit is the functor induced by the domain functors $C/x\to C$, while the comultiplication $\sum_x C/x\to \sum_{f\colon w\to x} C/w$ sends the $x$-component to the $1_x$-component via the identity functor $C/x \rightarrow C/x$. Dually, the left decalage comonad sends a category $C$ to the coproduct $\sum_x x/C$, with similar descriptions available for the counit and comultiplication.

Since the categories $C/x$ and $x/C$ are connected, a functor $\sum_x C/x\to \sum_x x/C$ is necessarily given by an assignment $c\mapsto tc$ on objects together with a functor $t\colon C/x\to tx/C$ for each $x$. Compatibility with the counit (on objects) means that the image under $t$ of an object $f\colon a\to x$ of $C/x$ should have the form $tf\colon tx\to a$. Functoriality, together with counit compatibility on morphisms means that if $fg=h$ then $g.th=tf$. Compatibility with the comultiplication requires a slightly more complicated calculation.

An object of $\Decr(C)$ has the form $f\colon a\to x$, and the comultiplication $\Decr(C)\to\Decr(\Decr(C))$ sends it to the composable pair $(1_x,f)$. Now $\Decr(t)\colon \Decr(\Decr(C))\to\Decr(\Decl(C))$ sends this to the composable pair  $(f,tf)$; which, as we have seen, must have composite $t1_x$. This composable pair can equally be seen as lying in $\Decl(\Decr(C))$, and finally applying $\Decl(t)$ gives the composable pair $(tf,t^21_x)$. Compatibility with comultiplication says that this should be equal to the composable pair $(tf,1_{tx})$, and this clearly says that $t^2(1_x)=1_{tx}$ for all objects $x$. We have only checked compatibility with the comultiplication on objects, but in fact no further condition is needed for compatibility on morphisms. We summarize this calculation as follows:

\begin{proposition}
To give duplicial structure to a small category $C$ is equivalently to give:
\begin{itemize}
\item for each object $x$ an object $tx$;
\item for each morphism $f\colon a\to x$ a morphism $tf\colon tx\to a$;
\end{itemize}
subject to the conditions that:
\begin{itemize}
\item $t^2(1_x)=1_{tx}$ for all objects $x$;
\item $g.t(fg)=tf$ for any composable pair $(f,g)$,
\end{itemize}
which we call the identity and functoriality conditions, respectively.
\end{proposition}

The next result gives a cleaner reformulation of these conditions. In its statement, recall that the inclusion $2$-functor $\mathbf{Gpd} \hookrightarrow \mathbf{Cat}$ has a left $2$-adjoint $\Pi_1$, whose counit at a small category $C$ is the functor $p \colon C \rightarrow \Pi_1(C)$ which freely adjoins an inverse for every arrow of $C$. The 2-dimensional aspect of the universal property means that, for any category $D$, the functor $[\Pi_1(C),D)\to [C,D]$ given by composition with $p$ is fully faithful.

\begin{theorem}
To give duplicial structure to a small category $C$ is equivalently to give a left adjoint in $\Cat$ for the functor $p\colon C\to \Pi_1(C)$.
\end{theorem}

\proof
First suppose that $p$ has a left adjoint $i\colon \Pi_1(C)\to C$ with counit $\epsilon\colon ip\to 1$ and unit $\eta\colon 1\to pi$; since $\Pi_1(C)$ is a groupoid, $\eta$ is invertible, and therefore $i$ is fully faithful. For each object $y\in C$ define $ty$ to be $ipy$, and for each morphism $f\colon x\to y$, define $tf\colon ipy\to x$ to be the composite
$$\xymatrix{
ipy \ar[r]^{i(pf)^{-1}} & ipx \ar[r]^{\epsilon_x} & x. }$$
Then $t(1_x)=\epsilon_x$ and so using the triangle identities twice yields
$$t^2(1_x)=t(\epsilon_x)=\epsilon_{tx}.i(p\epsilon_x)^{-1}=\epsilon_{ipx}.i\eta_{px} = 1_{ipx}$$
while for a composable pair $(g,f)$ we have
$$f.t(gf) = f.\epsilon_x.i(p(gf)^{-1}) = \epsilon_y.ip(f).i(pf)^{-1}.i(pg)^{-1} = \epsilon_y.i(pg)^{-1}=t(g)\rlap{ ;}$$
so this defines duplicial structure on $C$.

Conversely, if $C$ is equipped with duplicial structure there is an induced functor $G\colon C\to C$ sending an object $x$ to $tx$ and a morphism $f\colon x\to y$ to $t^2f\colon tx\to ty$. (This functor $G$ is the ``curious natural transformation'' of \cite{DwyerKan-cyclic} in another guise.) For each $x\in C$, write $\epsilon_x$ for the morphism $t(1_x)\colon tx\to x$.
Now $f.tf=t(1_y)$ by the functoriality condition, since $1_y f=f$; and replacing $f$ by $tf$ we also have $tf.t^2f=t(1_x)$. Combining these, $\epsilon_y.Gf=t1_y.t^2f=f.tf.t^2f=f.t1_x=f.\epsilon_x$ and so the $\epsilon_x$ are indeed natural. Furthermore, $G\epsilon_x=t^2(\epsilon_x)=t^3(1_x)=t(1_{tx})=\epsilon_{Gx}$ and so $(G,\epsilon)$ is a well-copointed endofunctor in the sense of~\cite{Kelly-transfinite}.

Next we show that for any $f\colon x\to y$, the morphism $Gf:=t^2f$ is invertible, with inverse $t(f.\epsilon_x)$. First observe that $\epsilon_x.t(f.\epsilon_x)=tf$ by the functoriality condition once again. Consequently we have
$$t(f.\epsilon_x).t^2(f)=t(f.\epsilon_x).t(\varepsilon_x.t(f.\varepsilon_x))=t(\epsilon_x)=t^2(1_x)=1_{tx}$$
using the functoriality condition again at the second step; this gives one of the inverse laws. By naturality of $\epsilon$ and the functoriality condition yet again we have
$$t^2f.t(f.\epsilon_x) = t^2f.t(\epsilon_y.t^2(f)) = t(\epsilon_y)=t^2(1_y)=1_{ty}$$
giving the other. Thus each $Gf$ is invertible. By the universal property of $\Pi_1(C)$, therefore, there is a unique functor $i\colon \Pi_1(C)\to C$ with $ip=G$. By the 2-dimensional aspect of the universal property of $\Pi_1(C)$, there is a unique natural transformation $\eta\colon 1\to pi$ with $\eta p\colon p\to pip$ equal to $(p\epsilon)^{-1}$, and so satisfying the triangle equation $p\epsilon.\eta p=1$. By the 2-dimensional aspect of the universal property once again, the other triangle equation $\epsilon i.i\eta=1$ will hold if and only if $\epsilon ip.i\eta p=1$ does, but by the calculation
$$\epsilon ip.i\eta p=\epsilon ip.(ip\epsilon)^{-1}=ip\epsilon.(ip\epsilon)^{-1}=1$$
this is indeed the case, and so $p$ does have a left adjoint.

It remains to show that these two processes are mutually inverse. First suppose that $C$ has duplicial structure $t$, and then construct a left adjoint $i\dashv p$ as above. The duplicial structure that this induces sends an object $x$ to $ipx=ix=tx$, and a morphism $f\colon x\to y$ to $\epsilon_x.i(pf)^{-1}$, where $i(pf)^{-1}=t(f.\epsilon_x)$. But now
$\epsilon_x.i(pf)^{-1} = \epsilon_x.t(f.\epsilon_x)=tf$ by the functoriality condition, and so we have recovered the original duplicial structure.

For the other direction, suppose first that $p$ has a left adjoint $i$ with counit $\epsilon$. Construct the induced duplicial structure $t$, and the left adjoint $i'$ and counit $\epsilon'$ induced by that. By the universal property of $\Pi_1(C)$ once again it will suffice to show that $ip=i'p$ and $\epsilon=\epsilon'$. For an object $x$, we have
$\epsilon'_x=t(1_x)=\epsilon_x.i(p1_x)^{-1}=\epsilon_x$, and so $\epsilon=\epsilon'$; this includes the fact that $ip$ and $i'p$ agree on objects, and so it remains only to show that they agree on morphisms. To see this, let $f\colon x\to y$ be a morphism, so that $i'pf\colon i'p x\to i'py$ is given by $t^2(f)\colon tx\to ty$. Now $tf=\epsilon_x.i(pf)^{-1}$, and so $ip(tf)^{-1}=ipipf.i(p\epsilon_x)^{-1}=ipipf.i\eta px=i\eta py.ipf$ and so finally
$i'pf=t^2f = \epsilon_{ipy}. i\eta py.ipf=ipf$.
\endproof

\begin{example}
If $C$ is a groupoid, then $p \colon C \rightarrow \Pi_1(C)$ is invertible, and so has a canonical left adjoint $p^{-1} \colon \Pi_1(C) \rightarrow C$. So every groupoid has a canonical duplicial structure.
\end{example}

\begin{example}\label{ex:groupoid-coreflection}
Suppose that there is a groupoid $G$ and a functor $i\colon G\to C$ with a right adjoint $r\colon C\to G$. By the universal property of $\Pi_1(C)$, there is a unique induced functor $q\colon \Pi_1(C)\to G$ with $qp=r$. By \cite[Proposition~1.3]{GabrielZisman}, this $q$ is an equivalence. Thus $p$ also has a left adjoint, and so $C$ has a duplicial structure.
\end{example}

\begin{remark}
We have seen that a category $C$ has duplicial structure just when $p\colon C\to \Pi_1(C)$ has a left adjoint. This will be paracyclic just when each $t_n$ is invertible, or equivalently just when each $t^{n+1}_n$ is invertible. Now the $t^{n+1}_n$ define the functor $ip\colon C\to C$; since $p$ is bijective on objects and $i$ is fully faithful, the composite $ip$ will be invertible if and only if $i$ and $p$ are both invertible, and this can happen only if $C$ is a groupoid.

For a groupoid, giving duplicial structure is equivalent to giving a left adjoint to the invertible $p \colon C \rightarrow \Pi_1(C)$; of course such a left adjoint is necessarily isomorphic to $p^{-1}$ and so in particular an equivalence. The duplicial structure will be paracyclic just when this left adjoint is in fact an invertible functor, and cyclic just when it is $p^{-1}$ as above. Thus, for a category $C$, the existence of paracyclic structure implies the existence of cyclic structure, but this does not mean that paracyclic structure on a category is necessarily cyclic. Furthermore, a groupoid can admit multiple cyclic structures, since there can be multiple choices of unit and counit for an adjunction $p^{-1}\dashv p$: in fact such choices correspond to choices of a natural isomorphism $1_G\cong 1_G$.
\end{remark}

\subsection{Duplicial structure on bicategories}

We next consider what it means to give duplicial structure on the nerve of a bicategory $B$ \cite{Street:categorical-structures}. Recall that this nerve is the simplicial set $NB$ in which:
\begin{itemize}
\item the 0-simplices are the objects of $B$;
\item the 1-simplices are the arrows $f\colon x\to y$ of $B$;
\item the 2-simplices are the 2-cells in $B$ of the form
$$\xymatrix{
& y \ar[dr]^{g} \dtwocell[0.6]{d}{\alpha} \\
x \ar[ur]^{f} \ar[rr]_{h} & & z ; }$$
\item the 3-simplices are the commuting diagrams of 2-cells of the form
$$\xymatrix{
(hg)f \ar[r]^{\cong} \ar[d]_{\alpha f} & h(gf) \ar[r]^{h\beta} & hk \ar[d]^{\gamma} \\
\ell f \ar[rr]_\delta && m }$$
in which the unnamed isomorphism is the relevant associativity constraint of $B$.
\end{itemize}
The face and degeneracy maps are as expected, and the higher simplices are determined by $3$-coskeletality. The assignment $B\mapsto NB$ is the object part of a fully faithful functor $N\colon \NLax\to\SSet$, where $\NLax$ is the category of bicategories and \emph{normal lax functors} between them---ones preserving identities on the nose, but binary composition only up to non-invertible $2$-cells $Fg.Ff \Rightarrow F(gf)$.  The first appearance in print we could find of the fact that this nerve functor is fully faithful was in \cite{BullejosFaroBlanco}. 

Once again, the decalage comonads on \SSet restrict to the full subcategory \NLax, and so it makes sense to speak of duplicial structure {\em on a bicategory}. Indeed the description of these restricted comonads is similar to the case of \Cat, except that rather than slice categories now we use ``lax slices''. For an object $x$ of a bicategory $B$, we write $B/x$ for the bicategory whose objects are morphisms $f\colon a\to x$ with codomain $x$, whose morphisms from $f\colon a\to x$ to $g\colon b\to x$ have the form
$$\xymatrix{
a \ar[rr]^{s} \ar[dr]_f^{~}="2"  && b \ar[dl]^{g}_{~}="1" \\
& x
\ar@{=>}"1";"2"^{\sigma} }$$
and whose 2-cells are defined in the evident way. Similarly the ``lax coslice'' $x/B$ has objects of the form $f\colon x\to a$, and morphisms from $f\colon x\to a$ to $g\colon x\to b$ of the form
$$\xymatrix{
& x \ar[dr]^{g}_{~}="2" \ar[dl]_{f}^{~}="1" \\
a \ar[rr]_s && b.
\ar@{=>}"1";"2"_{\sigma} }$$
We now define $\Decr(B)=\sum_x B/x$ and $\Decl(B)=\sum_x x/B$, with the actions on normal lax functors, and the counits and comultiplications given by a straightforward generalization of the corresponding definitions for \Cat.

Before giving our characterisation result, let us recall that a 2-cell in a bicategory as on the left in
$$\xymatrix{
& a \ar[dr]^{g}_{~}="2" \ar[dl]_{f}^{~}="1" \\
b \ar[rr]_h & & c
\ar@{=>}"1";"2"^{\alpha} } \qquad  \qquad
\xymatrix{
& a \ar[dr]^{g}_{~}="2" \ar[dl]_{k}^{~}="1" \\
b \ar[rr]_h & & c
\ar@{=>}"1";"2"^{\beta} }$$
is said to \emph{exhibit $f$ as a right lifting of $g$ through $h$} \cite{yoneda} if every 2-cell as right above factors as $\alpha . h\bar \beta$ for a unique 2-cell $\bar \beta \colon k \Rightarrow f$.

\begin{theorem}\label{thm:duplicial-bicategory}
 To equip a bicategory $B$ with duplicial structure is equivalently to give:
 \begin{enumerate}[(a)]
 \item for each object $x\in B$ an object $tx\in B$ and a morphism $\epsilon_x\colon tx\to x$;
\item for each morphism $f\colon a\to x$ in $B$ a morphism $tf\colon tx\to a$ and a 2-cell
$$\xymatrix{
& tx \ar[dr]^{\varepsilon_x}_{~}="2" \ar[dl]_{tf}^{~}="1" \\
a \ar[rr]_f & & x
\ar@{=>}"1";"2"^{\epsilon_f} }$$
exhibiting $tf$ as a right lifting of $\epsilon_x$ through $f$;
 \end{enumerate}
all subject to the conditions that:
\begin{enumerate}[(a)]
\addtocounter{enumi}{2}
\item $t1_x=\epsilon_x$;
\item $t^21_x=1_{tx}$;
\item $\epsilon_{t1_x} = 1_{t(1_x)}$.
\end{enumerate}
\end{theorem}

\proof
By redefining the composition with identity $1$-cells, any bicategory may be made isomorphic in \NLax to one in which identities are strict. Thus without loss of generality we may suppose that $B$ has strict identities.

Duplicial structure consists of a normal lax functor $t\colon \Decr(B)\to\Decl(B)$ which is compatible with the counit and comultiplication maps. As in the case of \Cat, since each $B/x$ and $x/B$ is connected, $t$ must be given by an assignment $x\mapsto tx$ on objects and normal lax functors $B/x\to tx/B$.

To give $t$ on objects compatibly with the counits is to give, for each $f\colon a\to x$, a morphism $tf\colon tx\to a$. To give $t$ on morphisms compatibly with the counits is to give, for each triangle as on the left below, a triangle as on the right.
$$
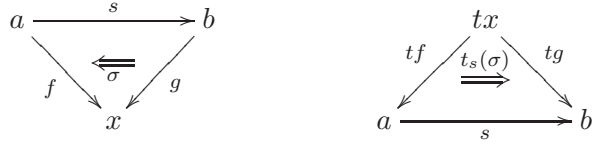
\xymatrix{
a \ar[rr]^s \ar[dr]_{f}^{~}="2" && b \ar[dl]^{g}_{~}="1" &&& tx \ar[dl]_{tf}^{~}="3" \ar[dr]^{tg}_{~}="4" \\
& x &&& a \ar[rr]_{s} && b
\ar@{=>}"1";"2"^{\sigma}
\ar@{=>}"3";"4"^{t_s(\sigma)}
}$$
The action of $t$ on 2-cells is unique if it exists, given the counit condition; it will exist just when, for all $\sigma\colon gs\to f$ and $\tau\colon s'\to s$, the diagram on the left commutes, where $\sigma'$ is defined as in the diagram on the right
$$  
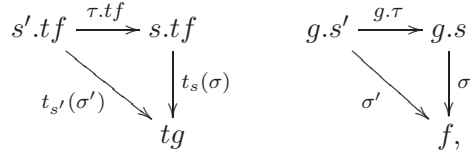
\xymatrix{
s'.tf \ar[r]^{\tau.tf} \ar[dr]_{t_{s'}(\sigma')} & s.tf \ar[d]^{t_s(\sigma)} \\
& tg } \quad\quad
\xymatrix{
g.s' \ar[r]^{g.\tau} \ar[dr]_{\sigma'} & g.s \ar[d]^{\sigma} \\ & f,}
$$
or, more compactly:
\begin{equation}
  \label{eq:sigma'}
t_{s'}(\sigma\circ(g.\tau)) = t_s(\sigma)\circ (\tau.tf)\rlap{ .}
\end{equation}
Since the components $\Decr(B)\to B$ and $\Decl(B)\to B$ of the counit are strict morphisms of bicategories, it follows that $t\colon \Decr(B)\to\Decl(B)$ will also be strict, which amounts to the requirements
\begin{equation}
  \label{eq:t-functorial}
  t_{1_a}(1_f) = 1_{tf} \quad \text{and} \quad t_{s'}(\sigma')\circ (s'.t_s(\sigma)) = t_{s's}(\sigma\circ\sigma's)
\end{equation}
for all $\sigma\colon gs\to f$ and $\sigma'\colon hs'\to g$.

It remains to see what the comultiplication axiom imposes. As in the case for \Cat, the only new condition appears at the level of objects of $\Decr(B)$, where it says that for any $f\colon a\to x$, we have
\begin{equation}
  \label{eq:delta-condition}
  t^21_x=1_{tx} \quad \text{and} \quad (t_{tf}(t_f(1_f))\colon tf.tt1_x\to tf)=1_{tf}\rlap{ .}
\end{equation}

So duplicial structure on a bicategory $B$ amounts to the assignments $x\mapsto tx$, $(f\colon a\to x)\mapsto (tf\colon tx\to a)$, and $(s,\sigma\colon gs\to f)\mapsto (t_s\sigma\colon s.tf\to tg)$; subject to the conditions expressed in \eqref{eq:sigma'}, \eqref{eq:t-functorial}, and \eqref{eq:delta-condition}. We now relate this to the structure in the statement of the theorem.

For any $x \in B$, we define $\varepsilon_x = t(1_x) \colon tx \rightarrow x$, and for any $f\colon a\to x$ in $B$, we define $\epsilon_f = t_f(1_f) \colon f.tf\to t1_x = \varepsilon_x$. Now in the conditions appearing in the theorem, (c) holds by construction, (d) holds by the first half of \eqref{eq:delta-condition}, while (e) holds by taking $f=1_x$ in the first half of \eqref{eq:t-functorial} and the second half of \eqref{eq:delta-condition}. Thus in order to show that a duplicial bicategory has all of the structure in the theorem, it remains only to show that $t_f(1_f)$ exhibits $tf$ as a right lifting of $t1_x$ through $f$; in other words, that for any $g\colon tx\to a$ and any $\phi\colon fg\to t1_x$, there is a unique $\psi\colon g\to tf$ which when pasted with $\epsilon_f$ gives $\phi$. But we may consider the pair $(g,\psi)$ as a morphism in $B/x$ from $t1_x$ to $f$, and so obtain  $t_g(\phi)\colon g.t^21_x\to tf$, and since $t^21_x=1_{tx}$, this gives our $\psi\colon g\to tf$. Pasting it with $\epsilon_f$ gives
\begin{align*}
\epsilon_f \circ f\psi &= t_f(1_f) \circ (f.t_g(\phi)) \\
  &= t_{fg} (\phi) \tag{by \eqref{eq:delta-condition}} \\
  &= t_{fg} ( 1_{t1_x}\circ (1_x.\phi)) \\
  &= t_{t1_x}(1_{t1_x})\circ \phi \tag{by \eqref{eq:sigma'}} \\
   &= \phi \tag{by (e)}
\end{align*}
which proves the existence of $\psi$. As for uniqueness, suppose that $\psi\colon g\to tf$ satisfies $\epsilon_f\circ f\psi=\phi$; that is, $t_f(1_f)\circ (f.\psi) = \phi$. Then
\begin{align*}
t_g(\phi) &= t_g(t_f(1_f)\circ (f.\psi)) \\
  &= t_{tf}(t_f(1_f))\circ (\psi.t^21_x) \tag{by \eqref{eq:sigma'}} \\
  &= \psi.t^21_x \tag{by \eqref{eq:delta-condition}} \\
  &= \psi \tag{by \eqref{eq:delta-condition}}
\end{align*}
giving uniqueness as required.

Thus a duplicial bicategory satisfies the conditions in the theorem.
For the converse, suppose that $B$ is equipped with structure as in the theorem; then we are given the assignments $x\mapsto tx$ and $(f\colon a\to x)\mapsto (tf\colon tx\to a)$, as well as the 2-cells $t_f(1_f)\colon f.tf\to \varepsilon_x$ satisfying the universal property of (b) as well as the conditions (c), (d), and (e). Given $\sigma\colon gs\to f$, if we are to have  \eqref{eq:t-functorial} and then \eqref{eq:sigma'} then
$$\epsilon_g\circ (g.t_s(\sigma))= t_g(1_g)\circ (g.t_s(\sigma)) = t_{gs}(\sigma) = \epsilon_f\circ (\sigma.tf)$$
and so $t_s(\sigma)$ is uniquely determined using the universal property of the right lifting 2-cell $\varepsilon_g$. It remains to check that if we define $t_s(\sigma)$ in this way, then \eqref{eq:sigma'}, \eqref{eq:t-functorial}, and \eqref{eq:delta-condition} do indeed hold.

Since $\epsilon_g\circ(g.t_s(\sigma))\circ(g.\tau.tf)=\epsilon_f\circ(\sigma.tf)\circ(g.\tau.tf)$, the composite $t_s(\sigma)\circ(\tau.tf)$ satisfies the defining property of $t_{s'}(\sigma\circ(g.\tau)$, and so \eqref{eq:sigma'} holds. Similarly
$\epsilon_h\circ(h.t_{s'}(\sigma'))\circ(h.s'.t_s(\sigma)) = \epsilon_g\circ(\sigma'.tg)\circ(h.s'.t_s(\sigma))=\epsilon_g\circ(g.t_s(\sigma))\circ(\sigma'.s.tf)=\epsilon_f\circ(\sigma.f)\circ(\sigma'.s.tf)$ and so $t_{s'}(\sigma')\circ(s'.t_s(\sigma))$ satisfies the defining property of $t_{s's}(\sigma\circ \sigma's)$; while $1_{tf}$ clearly satisfies the defining property of $t_{1_a}(1_f)$, and so \eqref{eq:t-functorial} holds.

The first half of \eqref{eq:delta-condition} is just (d); as for the second half, it says that $t_{tf}(\epsilon_f)=1_{tf}$, and the defining property of $t_{tf}(\epsilon_f)$ is that
$\epsilon_f\circ(f.t_{tf}(\epsilon_f))=\epsilon_{t1_x}\circ(\epsilon_f.t^21_x)$; but $t^21_x=1_{tx}$ by (d), and $\epsilon_{t1_x}=1_{t1_x}$ by (e), thus the right hand side becomes $\epsilon_f$, and clearly $\epsilon_f\circ 1_{tf}=\epsilon_f$, whence the result.
\endproof

\subsection{Monoidal categories}

A monoidal category can be thought of as a one-object bicategory, and as such it has a nerve: there is a unique 0-simplex, the 1-simplices are the objects of the monoidal category, the 2-simplices consist of three objects $X,Y,Z$ and a morphism $f\colon X\ox Y\to Z$, and so on. Thus the monoidal categories determine a full subcategory of \SSet, with the morphisms being the (lax) monoidal functors which are strict with respect to the unit. It is not the case that the decalage comonads restrict to this full subcategory: the decalage of a one-object bicategory will generally have many objects, indeed an object of the decalage will be a morphism of the monoidal category. Nonetheless, we can ask what it is to have duplicial structure on a monoidal category, thought of as a one-object bicategory.

Reading off directly from Theorem~\ref{thm:duplicial-bicategory}, we see that, for a monoidal category $C$ with tensor product $\ox$ and unit $i$, duplicial structure on $C$ consists of the following:
\begin{enumerate}[(a)]
\item an object $d$ (corresponding to $\epsilon_x$ for the unique object $x$ of the bicategory);
\item for each object $x$, a right internal hom $[x,d]$, by which we mean an object equipped with a morphism $\epsilon_x\colon x\ox [x,d]\to d$ inducing a bijection $C(x\ox-,d)\cong C(-,[x,d])$
\end{enumerate}
subject to conditions which we now enumerate. First of all, we require that the internal hom $[i,d]$ be $d$ itself. This is not a restriction in practice, since in any monoidal category and any object $x$ the internal hom $[i,x]$ exists and may be taken to be $x$. The more serious requirement is that the (chosen) hom $[d,d]$ is $i$, with counit $d\ox i\to d$ given by the unit isomorphism of the monoidal category. In fact the real condition here is that the map $i\to [d,d]$ induced by the unit isomorphism $d\ox i\to d$ is invertible; when this is the case we may always redefine $[d,d]$ as required.

One formulation of the notion of $*$-autonomous category \cite[Definition~2.3]{Barr-NonsymmetricStarAutonomous} is a monoidal category $C$ equipped with an equivalence $(-)^*\colon C\to C\op$ and a natural isomorphism $C(x,y^*)\cong C(i,(x\ox y)^*)$, with $i$ the unit. Using the natural isomorphism, we may construct further isomorphisms $C(x,y^*)\cong C(i,(x\ox y)^*)\cong C(i,(x\ox y\ox i)^*)\cong C(x\ox y,i^*)$, and so $y^*$ must in fact be given by $[y,i^*]$. Conversely, if $C$ is a monoidal category with all internal homs $[x,d]$ for a given object $d$, then there is a functor $(-)^*\colon C\to C\op$ sending $x$ to $[x,d]$, and a natural isomorphism $C(x,y^*)\cong C(i,(x\ox y)^*)$.

Thus both duplicial monoidal categories and $*$-autonomous categories involve a monoidal category $C$ and an object $d$ for which the right internal homs $[x,d]$ exist. The difference is that $*$-autonomous categories require the functor $[-,d]$ to be an equivalence, while duplicial monoidal categories require the canonical map $i\to [d,d]$ to be invertible. But in fact for a $*$-autonomous category the canonical map $i\to [d,d]$ is always invertible \cite[Section~6]{Barr-NonsymmetricStarAutonomous} and so any $*$-autonomous category has duplicial structure.

On the other hand, if the monoidal category $C$ has paracyclic structure, then the functor $[-,d]\colon C\to C\op$ will be (not just an equivalence but) invertible, and so $C$ will be $*$-autonomous.

The monoidal category $C$ will have cyclic structure when applying $[-,d]$ twice gives the identity; a mild generalization of this is that $[-,d]$ also gives {\em left} internal homs, in which case the $*$-autonomous category $C$ is said to have {\em cyclic dualizing object} \cite[Section~4]{Barr-NonsymmetricStarAutonomous}.



\bibliographystyle{plain}


\begin{thebibliography}{10}

\bibitem{Barr-NonsymmetricStarAutonomous}
Michael Barr.
\newblock Nonsymmetric {${}^\ast$}-autonomous categories.
\newblock {\em Theoret. Comput. Sci.}, 139(1-2):115--130, 1995.

\bibitem{Beck:dist-law}
Jon Beck.
\newblock Distributive laws.
\newblock In {\em Sem. on Triples and Categorical Homology Theory (ETH,
  Z\"urich, 1966/67)}, pages 119--140. Springer, Berlin, 1969.

\bibitem{BohmStefan1}
Gabriella B{\"o}hm and Drago{\c{s}} {\c{S}}tefan.
\newblock ({C}o)cyclic (co)homology of bialgebroids: an approach via
  (co)monads.
\newblock {\em Comm. Math. Phys.}, 282(1):239--286, 2008.

\bibitem{BohmStefan2}
Gabriella B{\"o}hm and Dragos {\c{S}}tefan.
\newblock Examples of para-cocyclic objects induced by {BD}-laws.
\newblock {\em Algebr. Represent. Theory}, 12(2-5):153--180, 2009.

\bibitem{BohmStefan3}
Gabriella B{\"o}hm and Drago{\c{s}} {\c{S}}tefan.
\newblock A categorical approach to cyclic duality.
\newblock {\em J. Noncommut. Geom.}, 6(3):481--538, 2012.

\bibitem{BullejosFaroBlanco}
M.~Bullejos, E.~Faro, and V.~Blanco.
\newblock A full and faithful nerve for 2-categories.
\newblock {\em Appl. Categ. Structures}, 13(3):223--233, 2005.

\bibitem{Burroni-NDA}
Elisabeth Burroni.
\newblock Alg\`ebres non d\'eterministiques et {$D$}-cat\'egories.
\newblock {\em Cahiers Topologie G\'eom. Diff\'erentielle}, 14(4):417--475,
  480--481, 1973.
\newblock Conf{\'e}rences du Colloque sur l'Alg{\`e}bre des Cat{\'e}gories
  (Amiens, 1973), I.

\bibitem{Connes-cyclic}
Alain Connes.
\newblock Cohomologie cyclique et foncteurs {${\rm Ext}\sp n$}.
\newblock {\em C. R. Acad. Sci. Paris S\'er. I Math.}, 296(23):953--958, 1983.

\bibitem{Day-convolution}
Brian Day.
\newblock On closed categories of functors.
\newblock In {\em Reports of the Midwest Category Seminar, IV}, Lecture Notes
  in Mathematics, Vol. 137, pages 1--38. Springer, Berlin, 1970.

\bibitem{LaxCentre}
Brian Day, Elango Panchadcharam, and Ross Street.
\newblock Lax braidings and the lax centre.
\newblock In {\em Hopf algebras and generalizations}, volume 441 of {\em
  Contemp. Math.}, pages 1--17. Amer. Math. Soc., Providence, RI, 2007.

\bibitem{DwyerKan-cyclic}
W.~G. Dwyer and D.~M. Kan.
\newblock Normalizing the cyclic modules of {C}onnes.
\newblock {\em Comment. Math. Helv.}, 60(4):582--600, 1985.

\bibitem{Elmendorf-CyclicDuality}
A.~D. Elmendorf.
\newblock A simple formula for cyclic duality.
\newblock {\em Proc. Amer. Math. Soc.}, 118(3):709--711, 1993.

\bibitem{GabrielZisman}
P.~Gabriel and M.~Zisman.
\newblock {\em Calculus of fractions and homotopy theory}.
\newblock Ergebnisse der Mathematik und ihrer Grenzgebiete, Band 35.
  Springer-Verlag New York, Inc., New York, 1967.

\bibitem{GetzlerJones-paracyclic}
Ezra Getzler and John D.~S. Jones.
\newblock The cyclic homology of crossed product algebras.
\newblock {\em J. Reine Angew. Math.}, 445:161--174, 1993.

\bibitem{Kelly-transfinite}
G.~M. Kelly.
\newblock A unified treatment of transfinite constructions for free algebras,
  free monoids, colimits, associated sheaves, and so on.
\newblock {\em Bull. Austral. Math. Soc.}, 22(1):1--83, 1980.

\bibitem{property}
G.~M. Kelly and Stephen Lack.
\newblock On property-like structures.
\newblock {\em Theory Appl. Categ.}, 3:No. 9, 213--250 (electronic), 1997.

\bibitem{Slevin2}
Niels Kowalzig, Ulrich Kraehmer, and Paul Slevin.
\newblock Cyclic homology arising from adjunctions.
\newblock {\em Theory and Applications of Categories}, 30:1067--1095, 2015.

\bibitem{Slevin1}
Ulrich Kraehmer and Paul Slevin.
\newblock Factorisations of distributive laws.
\newblock  {\em J. Pure Appl. Algebra}, in press, 2015.

\bibitem{codescent}
Stephen Lack.
\newblock Codescent objects and coherence.
\newblock {\em J. Pure Appl. Algebra}, 175(1-3):223--241, 2002.

\bibitem{skew}
Stephen Lack and Ross Street.
\newblock Skew monoidales, skew warpings and quantum categories.
\newblock {\em Theory Appl. Categ.}, 26:385--402, 2012.

\bibitem{Loday-cyclic}
Jean-Louis Loday.
\newblock {\em Cyclic homology}, volume 301 of {\em Grundlehren der
  mathematischen Wissenschaften}.
\newblock Springer, Berlin Heidelberg, 1992.

\bibitem{CWM}
Saunders MacLane.
\newblock {\em Categories for the working mathematician}.
\newblock Springer-Verlag, New York, 1971.

\bibitem{Schauenburg-DualsAndDoubles}
Peter Schauenburg.
\newblock Duals and doubles of quantum groupoids ({$\times\sb R$}-{H}opf
  algebras).
\newblock In {\em New trends in {H}opf algebra theory ({L}a {F}alda, 1999)},
  volume 267 of {\em Contemp. Math.}, pages 273--299. Amer. Math. Soc.,
  Providence, RI, 2000.

\bibitem{ftm}
Ross Street.
\newblock The formal theory of monads.
\newblock {\em J. Pure Appl. Algebra}, 2(2):149--168, 1972.

\bibitem{FibBicCorr}
Ross Street.
\newblock Correction to: ``{F}ibrations in bicategories'' [{C}ahiers
  {T}opologie {G}\'eom.\ {D}iff\'erentielle {\bf 21} (1980), no.\ 2, 111--160;
  {MR}0574662 (81f:18028)].
\newblock {\em Cahiers Topologie G\'eom. Diff\'erentielle Cat\'eg.},
  28(1):53--56, 1987.

\bibitem{Street:categorical-structures}
Ross Street.
\newblock Categorical structures.
\newblock In {\em Handbook of algebra, Vol.\ 1}, pages 529--577. North-Holland,
  Amsterdam, 1996.

\bibitem{yoneda}
Ross Street and Robert Walters.
\newblock Yoneda structures on 2-categories.
\newblock {\em J. Algebra}, 50(2):350--379, 1978.

\bibitem{Szlachanyi-skew}
Korn{\'e}l Szlach{\'a}nyi.
\newblock Skew-monoidal categories and bialgebroids.
\newblock {\em Adv. Math.}, 231(3-4):1694--1730, 2012.

\bibitem{Verity-complicial}
Dominic Verity.
\newblock Complicial sets characterising the simplicial nerves of strict
  {$\omega$}-categories.
\newblock {\em Mem. Amer. Math. Soc.}, 193(905):xvi+184, 2008.

\end{thebibliography}



\end{document}